# Symmetry Classification and Enumeration of Square-Tile Sikku Kolams

Venkatraman Gopalan

*Department of Materials Science and Engineering, The Pennsylvania State University, University Park, PA, 16802*

vgopalanpsu@gmail.com

**Abstract**

Pulli kolam is an ancient mathematical art form that is still practiced today in south India by over a quarter million people. 'Pulli' in the Tamizh language means 'dots'. A specific type of pulli kolam is sikku kolam where a series of dots are placed using rice flour and lines are drawn around them following three simple rules. Here we focus on sikku kolams assembled from square tiles. Kolams are classified by their point group symmetry and formulae are derived for counting the number of kolams of each symmetry type. Some interesting constraint problems are posed that could help formulate new games based on this elegant art form.

## 1 Introduction

Pulli kolam has been practiced for thousands of years in southern India, presently the states of Tamil Nadu, Andhra Pradesh, Telangana, Kerela and Karnataka, with a combined current population of over 260 million people. A particular type of pulli kolam is called sikku kolam (which translates in Tamizh to 'knot kolam') or kambi kolam (which translates to 'wire kolam'). This article focuses on the sikku kolams constructed from square tiles. We will refer to them as square-tile kolams. However, note that although the tiles themselves are squares, the assembled kolam need not be square; it can be any shape you are able to make by assembling square tiles following these three rules: (1) *All dots are individually encircled*, (2) *line crossings are discrete points, and* (3) *there are no loose ends in the lines drawn.*

**2 Six basic kolam tiles**

The six types of square tiles considered in this article are shown in Figure 1. They have been given descriptive names as indicated. Each tile has a single dot at its center with segments of one or more curves encircling it. Loose ends, in the case of these six tiles, occur only at midpoints of edges of the tiles and are indicated by grey arrows. The Circle, Drop, Eye, Door, Fan, and Diamond have respectively, 0, 1, 2, 2, 3, and 4 loose ends. When loose ends of two tiles come together, they annihilate each other, and the corresponding lines cross continuously across the two tiles. Note that there are many more tile types beyond these six basic types, indeed infinitely many. We will cover such an extension of the tile types in future works.

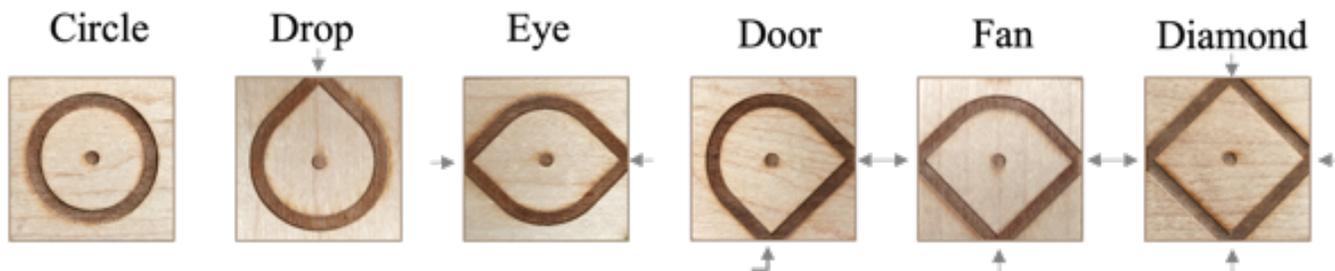

Figure 1. Six types of square kolam tiles considered in this article. The grey arrows indicate loose ends. These tiles were designed by the author in collaboration with Dr. Amy Alznauer, and were physically fabricated by Christopher Danielson.

The tiles in Figure 1 are designed in such a way that two of the three rules for a kolam are automatically satisfied, namely, (1) all the dots are individually encircled and (2) lines cross only at discrete points (the midpoints of edges of the tiles). For each square tile, there are exactly six possibilities for the numbers of loose ends it has: 0, 1, 2 (on adjacent sides), 2 (on opposite sides), 3, or 4, and Figure 1 has exactly one tile of each type. Therefore in assembling these tiles, there is only one remaining rule to pay attention to: (3) *no loose ends*. This is done by arranging the tiles to get rid of all loose ends by pairing loose ends that annihilate each other. Any assembly of these tiles that has no loose ends could be considered a valid pulli kolam.

**3 Prior literature**

Many previous pioneering works over the past five decades have provided mathematical insights into the form of a kolam. These include converting kolams into numbers and linear diagrams (Yanagisawa & Nagata, 2007), using graph, picture, and array grammers (Siromoney et al., 1974; Siromoney & Siromoney 1987; Subramanian & Siromoney, 1987; David et al. 2007; Pradella et al., 2011; Siromoney et. al., 1976; Siromoney et. al., 1983; Demaine et al., 2007; Subramanian, 1979), extended pasting schemes (Robinson, 2007), morphism of monoids (Allouche et al., 2006), L- and P-systems (Prusinkiewicz & Hanan, 1989; Subramanian, 2007), gestural lexicons (Waring, 2012), knot theory (Ishimoto, 2007), and mirror curves (Gerdes, 1999). Yanagisawa & Nagata begin with five rules for a kolam (instead of the three rules here) and define square unit tiles that can be assembled into larger kolams. They define two types of nearest-neighbor interactions between dots (line crossing, 1, or uncrossing, 0) and convert these tiles into binary number arrays. In Gopalan & VanLeeuwen, 2015, a purely topological approach is taken, specifying only three mandatory rules for defining a kolam. They generalize the ideas to any arbitrary arrangement of dots (not necessarily square arrays), generalize 'interactions' between any two dots (instead of only the nearest or next nearest neighbors), and to three or more 'bonds' between an interacting pair of dots. Their work suggests that for a given number of dots, $N$, there are a limited number of parent kolam types from which all other kolams originate. All parent kolams within a parent kolam type are homotopic (or topologically equivalent).

The present work is closely related to the approach of Gopalan & VanLeeuwen, where we restrict the type of kolam studied to those create using square tiles. The interactions between dots is restricted to only the nearest neighbor dots. This was also the specific case studied by Yanagisawa & Nagata, 2007 and Nagata, 2015; however they place more than one dot per tile and discuss the kolams in terms of linear diagrams and knot theory. In contrast, our work classifies the square tile kolams in terms of group theory, assigning the point groups well known in crystallography. The enumeration of different kolam types also uses these symmetry classifications. An example of such symmetry classification and enumeration is in the study of rosettes of Truchet tiles by Hall, Almeida & Teixeira, 2019.

**4 Eight symmetry types of square-tile kolams**

Using the six basic kolam tiles shown in Figure 1, we can create kolams. All such kolams created from square tiles can be classified into one of the eight point group symmetries depicted in Figure 2.

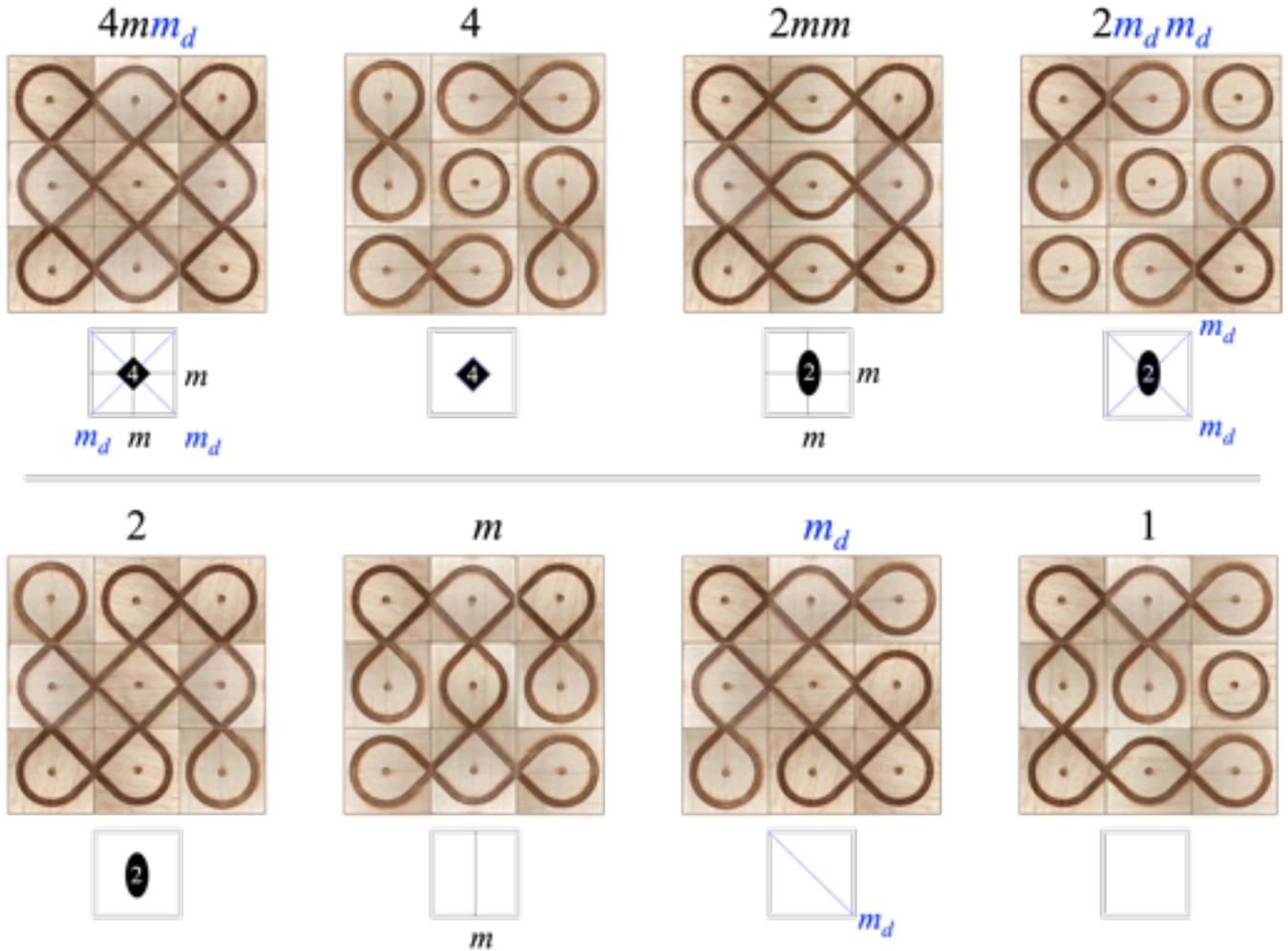

Figure 2. Symmetries of all square-tile kolams can be classified into one of the eight point groups illustrated above. The grey square next to each kolam represents the entire kolam and shows the kolam's symmetry elements. Lines depict mirrors, black squares depict 4-fold rotation centers, and lozenges depict 2-fold rotation centers. Each kolam is created from a subset of the six basic tile types shown in Figure 1.

Groups are mathematical sets with specific required properties (Washburn & Crowe, 1991). For a given kolam, the symmetry group is the set of all possible symmetry operations one can perform on the kolam that will leave the kolam unchanged. A point group specifies that all the symmetry operations must leave at least one point in the kolam unmoved. In our case, it will typically be the center of the kolam about which one can rotate the kolam, and through which all mirrors pass. Rotations and mirror reflections are the point group symmetry operations that are relevant to the 2-dimensional kolams. Here is some helpful terminology:

(1) A *p*-fold rotation axis means that the kolam can be rotated about this axis by 360°/*p* and the kolam will still look the same, where *p* is a natural number. In the square-tile kolams, 4-fold (90°), 2-fold (180°), and 1-fold (360°) rotation axes are the only ones possible.

(2) A mirror plane is represented by *m* and will reflect the coordinate perpendicular to the mirror while leaving unchanged the other coordinates parallel to the mirror. All mirrors of symmetry must pass through the center of a kolam. In the square-tile kolams, there are two types of mirrors: the ones parallel to the edges of the kolam (denoted in Figure 2 by black lines labeled *m*) and the ones lying on the diagonals of the kolam (denoted in Figure 2 by $m_d$ lines).

In this work, when we speak of the symmetry of a kolam, we include not only the dots and the curves, but the square tiles they are imprinted on as well. The largest possible symmetry point group for a square-tile kolam has crystallographic notation 4*mm*$_d$, as seen in Figure 2. The grey box next to the kolam titled *4mm*$_d$ schematically represents the entire kolam, and the various symmetry operations are symbolically represented within this grey box. The 4-fold rotation axis is perpendicular to the plane of the paper at the center of this kolam and allows a rotation of the entire kolam by 360°/4 = 90°, after which the kolam looks identical to its state before the rotation. The presence of a 4-fold rotation symmetry automatically implies the presence of a 2-fold rotation axis as well. All kolams also trivially have 1-fold rotation symmetry. There are two sets of mirrors also passing through the center of this kolam: the ones indicated in black (*m* are parallel to the edges of the kolam), while the ones in blue (namely, $m_d$) lie on the diagonals. The two black mirrors are symmetry-equivalent (conjugate) since the 4-fold axis can rotate one of them into the other. Similarly, the two blue mirrors are equivalent for the same reason. However, the black and the blue mirrors are not equivalent since no symmetry operation within the point group can transform a black mirror into a blue one or vice versa. In group theory, the sets of mirrors *m* and $m_d$ are said to belong to different *classes*. Note that the presence of any two of these symmetry operations (4-fold rotation, reflection in *m*, and reflection in $m_d$) would automatically generate the third.

With that detailed explanation, the other seven symmetry point groups depicted in Figure 2 should be clear, namely: 4, 2*mm*, 2$m_d m_d$, 2, *m*, $m_d$, and 1. Note that the point group 2*mm* has two inequivalent mirrors (i.e. they belong to different classes) since no symmetry operation within the group can transform one of those mirrors into the other one. The point group 1 has no symmetry operations other than the 1-fold, namely a rotation by 360° about an axis located at its center; the 1-fold axis is not labeled. As an aside, the pair of point groups, $m$ and $m_d$ would not be separately distinguished in crystallography, but we choose to distinguish them here in kolams for practical and artistic reasons. The same is true for the pair $2mm$ and $2m_d m_d$.

**5 Kolam count by symmetry group**

Of the infinitely many possible arrangements of the tiles in Figure 1, we consider just two basic types of square-tile kolam templates for enumeration. These are (1) single rectangle (1R) with $kl$ tiles (and hence $kl$ dots), and (2) two rectangle (2R), with $kl + (k-1)(l-1)$ tiles (and hence $kl + (k-1)(l-1)$ dots). Two examples are shown in Figure 3; these have symmetry type $2mm$ (see center schematic). However, all other symmetry types in Figure 2 are also possible. For a square template, $k = l$.

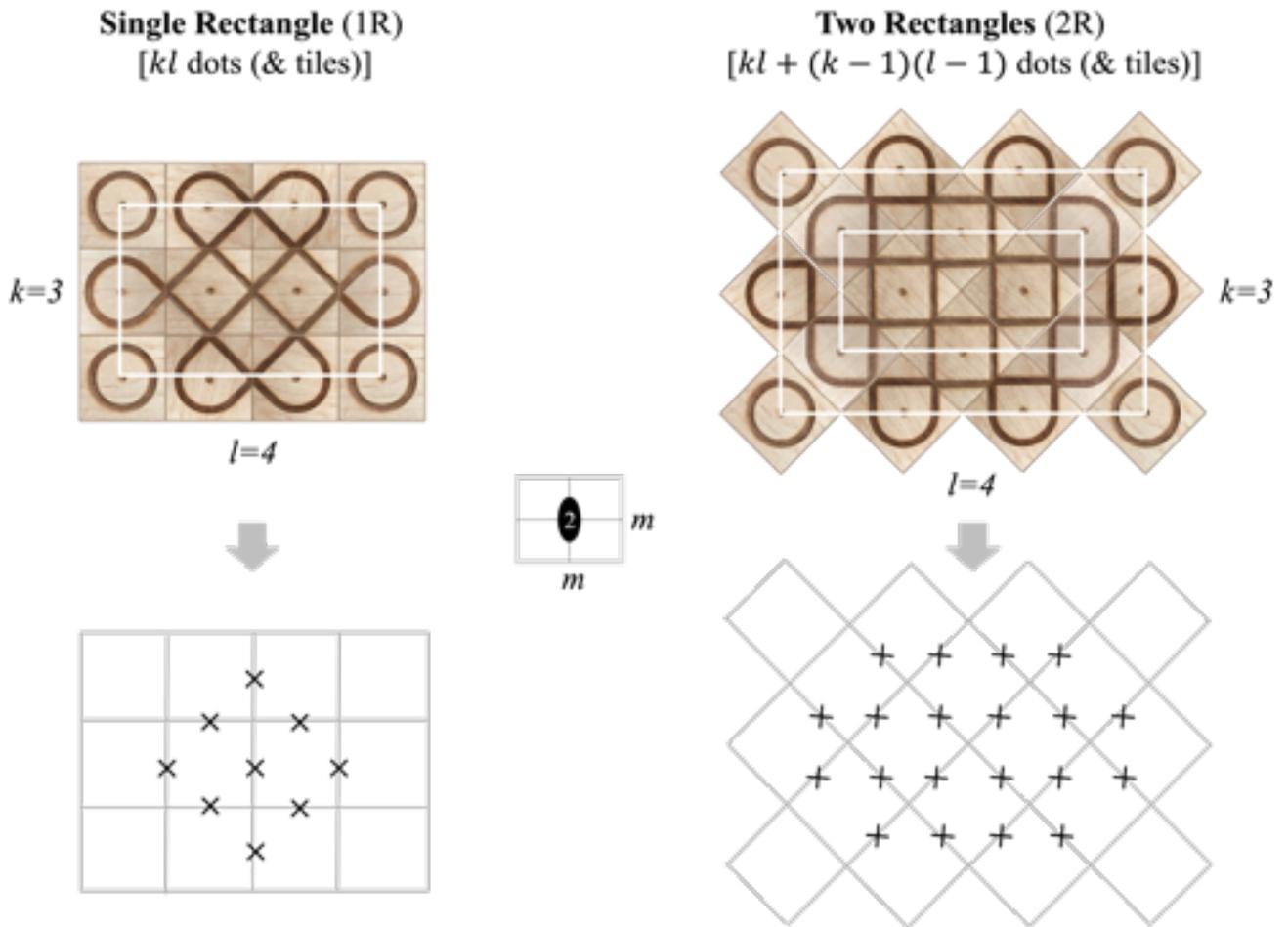

Figure 3. Two basic kolam templates considered for enumeration. The single rectangle (1R) template (images on the left) has square tiles with edges aligned horizontally and vertically and has $kl = 3 \times 4 = 12$ dots (its 3 x 4 rectangle outlined in white). The two rectangles (2R) template (images on the right) has square tiles with edges aligned diagonally and has $kl + (k-1) \times (l-1) = 3 \times 4 + 2 \times 3 = 18$ dots (its $3 \times 4$ rectangle and $2 \times 3$ rectangle outlined in white). The bottom two depictions outline in grey the tiles in the kolams above them and × symbols indicate where lines cross in the kolams.

Since there is exactly one tile for each of the possible numbers and positions of line crossings for each square in the grey grid, all of the information needed to construct the two kolams at the top of Figure 3 is contained in the bottom depiction of line crossings. One way to enumerate the possible kolams is therefore to count the total number of shared edges and assign to each shared edge either a line crossing, '×' or no crossing. Without any symmetry considerations, (i.e. for the point group 1 in Figure 2), each shared edge can be independently specified. For other symmetry groups, a subset of these can be independently specified to construct the whole kolam.

Thus, depending on the desired symmetry of a kolam, only a minimum number of shared edges, $E_s^{PG}$, need to be marked with one of two states, an × or "no ×", to produce the whole kolam. (Here '$PG$' labels the point group of the kolams being enumerated). For a kolam with no symmetry (symmetry group of 1), the state of every shared edge would need to be specified. But for all others, the state of only some of the shared edges would need to be specified, since the symmetry operation(s) would then specify the state of the unspecified shared edges. For example, for a kolam of point group symmetry $4mm_d$ with the 1R template and $n^2$ tiles, there are $2n(n-1)$ shared edges in total; of these, the state of only $n^2/4$ shared edges need to be specified if $n$ is even, and the state of only $(n^2-1)/4$ shared edges need to be specified if $n$ is odd. Thus with $n = 5$, $E_s^{4mm_d} = 6$. It is not difficult to find similar general expressions for the 1R and 2R templates shown in Figure 3 for each of the 8 point groups illustrated in Figure 2; these are given in Table 1 in the Appendix. One can now calculate the number of *unique* kolams for every symmetry group as follows:

$$N(PG) = 2^{E_s^{PG}} \qquad (1)$$

Here, the superscript $PG$ can be one of the eight point group symmetries shown in Figure 2. To get a sense for how many kolams are generated, Table 2 in the Appendix lists the number of unique kolams for some small kolam sizes for each template in Figure 3. Since the specification of $E_s^{PG}$ already accounts for the symmetry, it removes any duplicate kolams from being counted in Eq. (1) above. By duplicate, we mean that a completed kolam that is simply rotated or mirrored does not count as a unique kolam; hence these are not counted. Even so, the number of possible kolams becomes very large very quickly. The higher the symmetry, the fewer the number of unique kolams, and the lower the symmetry, the greater number of unique kolams. Figures 4 and 5 depict all the square-tile kolams with the largest symmetry, $4mm_d$ up to size 3 × 3 for the 2R templates and up to size 4 × 4 for the 1R templates. Other kolams can be similarly constructed by the technique shown in Figure 3.

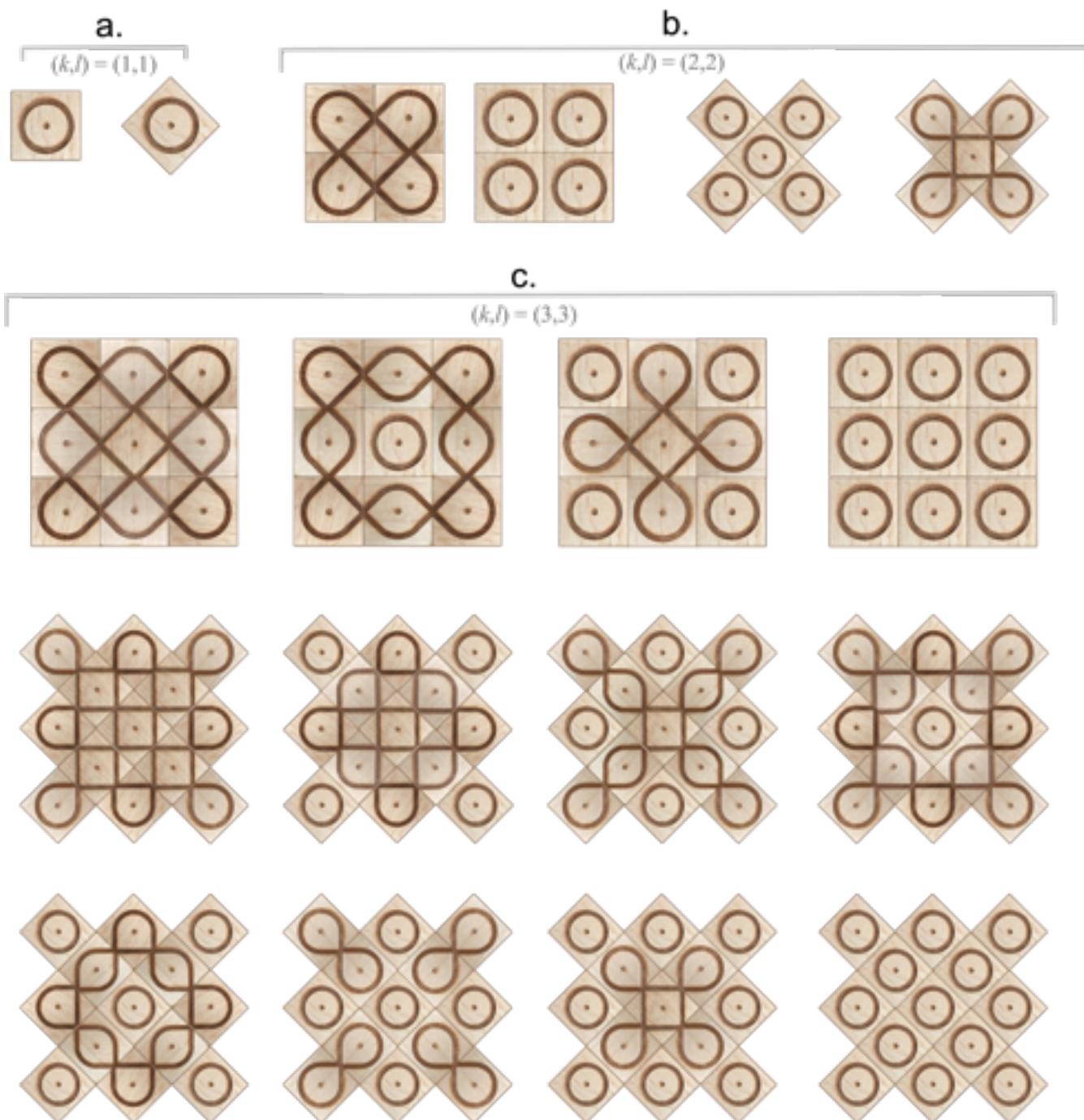

Figure 4. All the 1R and 2R template kolams with $4mm_d$ symmetry for (a) $k \times l = 1 \times 1$, (b) $k \times l = 2 \times 2$, and (c) $k \times l = 3 \times 3$. The counts are consistent with Table 2.

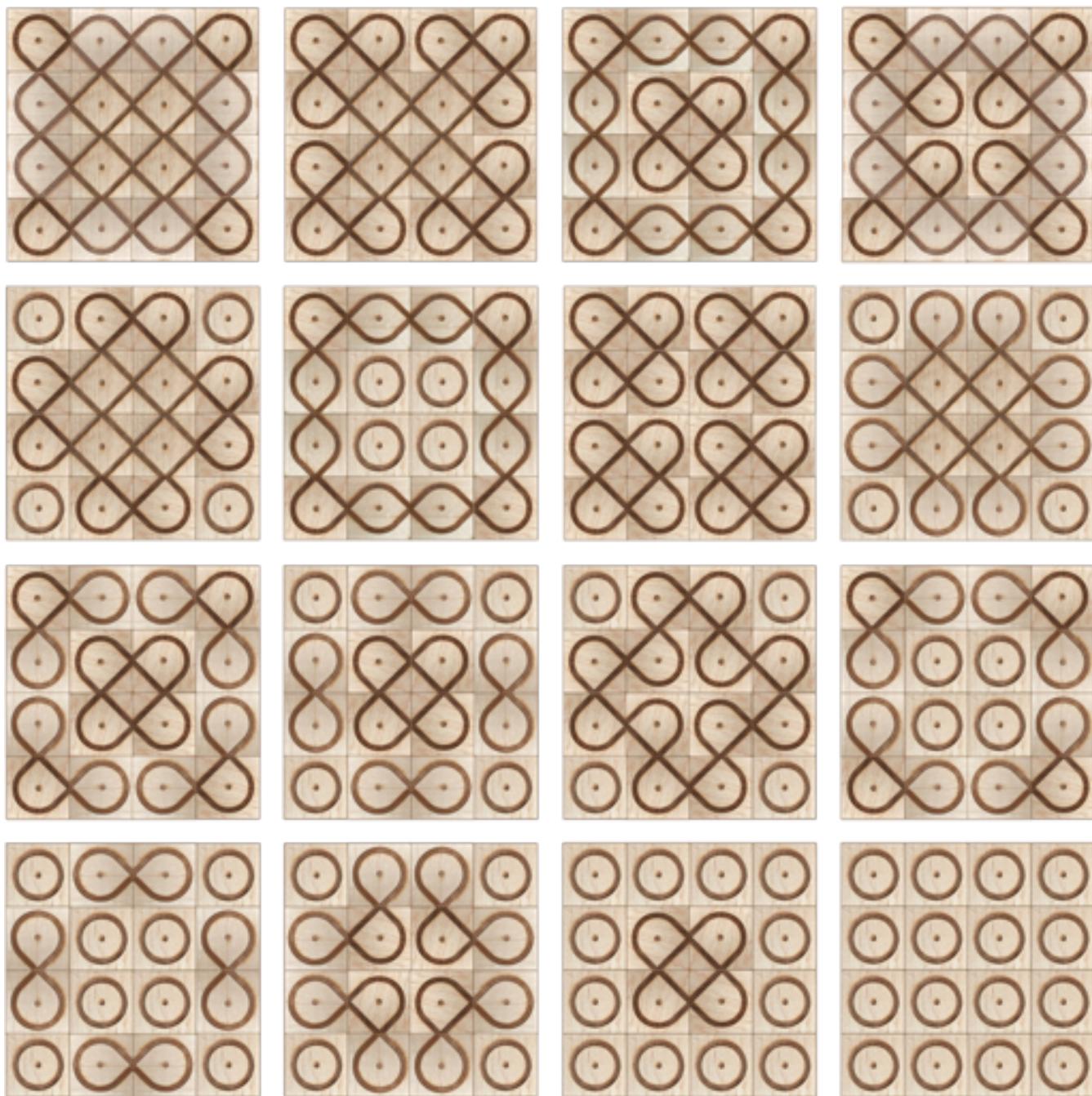

Figure 5. All the 1R template kolams with $4mm_d$ symmetry for $k \times l = 4 \times 4$. The count of 16 is consistent with Table 2.

## 6 Necessary and sufficient conditions for a kolam to exist

A sufficient condition for composing a kolam would be the three rules given in the Introduction. Two of these are automatically satisfied by the construction of the kolam tiles in Figure 1. Thus for the square-tile kolams, a sufficient condition is simply no loose ends.

We can now identify some necessary conditions for achieving no loose ends. We first show that

*A necessary condition to have no loose ends is:*
*the total number of Drop tiles and Fan tiles (in Figure 1) in the kolam must be even.* (2)

This condition makes it easy to determine when a kolam cannot be composed, given a specified collection of tiles. For example, suppose we are given the following tiles and asked to compose a kolam: one Circle, three Drops, two Eyes, two Fans, four Doors and one Square. Even before starting, we can tell that no kolam is possible with these tiles without leaving a loose end. This is because we have three drops and two fans, for a total of five, which is an odd number; hence no kolam is possible.

The proof for the necessary condition (2) is straightforward. Let $N_i$ represent the number of tiles of type $i$, where $i$ can be $C$ (Circle), $Dr$ (Drop), $E$ (Eye), $Do$ (Door), $F$ (Fan), or $Di$ (Diamond). Recall that these tiles respectively have 0, 1, 2, 2, 3, and 4 loose ends. Since the loose ends must occur in pairs to completely annihilate each other, we can conclude that:

$$0N_C + 1N_{Dr} + 2(N_E + N_{Do}) + 3N_f + 4N_{Di} = 2n, \qquad (n, \text{ a natural number}) \qquad (3)$$

This implies that $N_{Dr} + 3N_f$ must be even, and so $N_{Dr} + N_f$ must be even.

Even if we have an even number of loose ends that we can in principle annihilate, there needs to be a sufficient number of shared edges between tiles that allow for such annihilation. If $E_s^t$ is the total number of shared edges in a kolam, and a pair of loose ends annihilate each other at a shared edge, then it follows that in equation (3),

$$n \leq E_s^t < 2T \qquad (4)$$

where $T$ is the fixed number of tiles specified for composing a kolam, that is, $T = N_C + N_{Dr} + N_E + N_{Do} + N_f + N_{Di}$. The first inequality asserts that the total number of pairs of loose ends should be less than or equal to the number of shared edges, $E_s^t$, which by the second inequality should be less than the twice the total number of tiles, $T$, in the kolam. Even before starting to assemble a kolam from a set of square tiles, we can test for the necessary conditions given in (2) and the first inequality in (4). The second inequality in (4) is always satisfied and is not meant to be a "test"; it can be derived by noting that no kolam can have all edges of every tile to be shared, since the peripheral edges of a kolam by definition are not shared with any other tiles. Since, $T$ square-tiles have $4T$ edges, and $2E_s^t$ of these are part of shared edges, therefore, $2E_s^t < 4T$, and hence $E_s^t < 2T$.

## 7 Conditions imposed by symmetry

Additional constraints are imposed on the total number of tiles $T$ needed to create a kolam with specific symmetry.

- **4-fold rotation:** If the desired kolam is to have a 4-fold rotation axis, then when $T$ is even, there must be four copies of $T/4$ tiles, hence $T$ must be divisible by 4. When $T$ is odd, then there must be four copies of $(T-1)/4$ tiles; therefore $(T-1)$ must be divisible by 4. In addition, for odd $T$, there must be one tile that is placed in the center of the kolam. The 4-fold axis would pass through this center tile. For $T=1$, only one square tile, namely the Circle, can form a kolam by itself and has 4-fold symmetry when the square shape of the tile is included; if only the circle is considered without the square tile template, it has an $\infty$-fold rotation axis. A 4-fold symmetry kolam is not possible for any $T$ that is not the square of a natural number, hence $T = 2$ or $T = 3$ e.g. are ruled out. For $T \geq 4$, assuming the necessary conditions on $T$ are met, a minimum of $(T-1)/4 + 1$ tiles will need to be specified for $T$ odd and $T/4$ tiles for $T$ even to construct the whole kolam.

- **2-fold rotation:** If the desired kolam is to have a 2-fold rotation axis, then when $T$ is even, there must be two copies of $T/2$ tiles (which follows from the fact that $T$ is divisible by 2). When $T$ is odd, there must be two copies of $(T-1)/2$ tiles. Further, for odd $T$, there must be one odd tile that is placed at the center of the kolam, where the 2-fold axis is located. For $T = 1$, only the Circle tile qualifies. For $T \geq 2$, a minimum of $(T-1)/2 + 1$ tiles will need to be specified for $T$ odd and $T/2$ tiles specified for $T$ even in order to construct the whole kolam.

- **One mirror, $m$:** For a 1R template, to have one mirror parallel to one of the edge directions, ($k$ or $l$) with $T = k \times l$ tiles, and $T$ odd, then $(T-k)$ and $(T-l)$ are divisible by 2. In this case, only $(T+k)/2$ tiles will need to be specified to create the whole kolam with a mirror parallel to the $k$ edge. Since $T$ is odd, $k$ must be odd, and so there will be tiles cut by the mirror $m$ of the kolam; these tiles must have mirror symmetry themselves with respect to $m$, which rules out the Door tile. (The last two statements for $k$ can be repeated with $k$ replaced by $l$.)

When $T$ is even, and assuming that $l$ is even and $k$ is odd, only $T/2$ tiles need to be specified if the mirror $m$ is parallel to the $k$ edge; further, no tiles are cut by $m$. For the mirror $m$ to be parallel to the $l$ edge, $(T+1)/2$ tiles will need to be specified for creating the whole kolam, and the tiles cut by $m$ cannot be the Door tile. The case when $k$ and $l$ are both even have no tiles cut by $m$, and only $T/2$ tiles need to be specified.

For a 2R template, the same analysis as for 1R can be performed except the two cases of $T = kl$ tiles and $(k-1)(l-1)$ additional tiles must be considered separately. Distinct from the 1R template, 2R kolams will always have tiles cut by $m$, and these tiles must have mirror symmetry about the diagonal of the individual square tiles. This restricts the tiles cut by the mirror $m$ to Circles, Diamonds or Doors.

- **Two mirrors, $m$:** If there are only two mirrors in a square-tile 1R kolam, they must be orthogonal (perpendicular). Further, at their intersection is a 2-fold rotation axis; in other words, the point group symmetry will be $2mm$ as shown in Figure 2. The constraints imposed by the 2-fold rotation and the mirror $m$ as described above hold for this 1R kolam as well. The total number of tiles for a 1R kolam is $T = kl$. The minimum number of tiles that need to be specified for a 1R kolam with $2mm$ symmetry is $(T+l)/4$ tiles if $k$ is odd and $l$ is even, $(T+k)/4$ tiles if $l$ is odd and $k$ is even, $(T+l+k)/4$ if $k$ and $l$ are both odd, and $T/4$ tiles if $k$ and $l$ are both even. The tiles on the mirror lines of a 1R kolam cannot be the Door tiles.

For a 2R kolam, $T = kl + (k-1)(l-1)$. For a 2R kolam with $2mm$ symmetry, the same considerations as a 1R kolam apply, but need to be applied separately to a 1R kolam with $kl$ tiles in addition to a 1R kolam with $(k-1)(l-1)$ tiles, and then summing them up for the total minimum number of tiles that need to be specified. By performing this procedure, the minimum number of tiles that need to be specified for a 2R kolam with $2mm$ symmetry is $(T+k+l-1)/4$ for the cases of $k$-even, $l$-odd, or for $k$-odd, $l$-even. For the case of $k$-odd, $l$-odd, the minimum number of tiles needed to be specified is $(T+k+l)/4$, while for $k$-even, $l$-even, it is $(T+k+l-2)/4$. Distinct from the 1R, 2R kolams will always have tiles lying on the mirror lines, and these tiles must have mirror symmetry about the diagonal of the individual square tiles. This restricts the tiles that are cut by the mirror $m$ in 2R kolams to be either Circles, Diamonds or Doors.

- **One diagonal mirror, $m_d$:** A 1R kolam with $T = kl$ tiles can have a diagonal mirror $m_d$ of symmetry only if $k = l$ and hence $T = k^2$ must be a square number. Further, only Circle, Door and Diamond tiles can be along the diagonals of such a 1R kolam. Also, $(T-k)$ must be an even number, and a minimum of $(T + k)/2$ tiles are needed to complete the full kolam.

For a 2R template with one diagonal mirror $m_d$ parallel to one of the diagonals of the template, $T = k^2 + (k − 1)^2$ is a sum of two square numbers. Further, Doors cannot be on $m_d$. A 2R kolam with one $m_d$ mirror requires $(T + 2k + 1)/2$ tiles to complete the entire kolam.

• **Two diagonal mirrors, $m_d$:** The point groups $2m_d m_d$ and $4mm_d$ each have a pair of diagonal mirrors that are orthogonal. These kolams need to have a square number of tiles: $T = k^2$ for 1R and $T = k^2 + (k − 1)^2$ for the 2R. For the $2m_d m_d$ and $4mm_d$ 1R kolams, only Circle, Door, and Diamond tiles can be along the diagonals. The minimum number of tiles that need to be specified to construct the whole $2m_d m_d$ 1R kolam is given by $(T + 2k + 1)/4$ if $k$ is odd, and $(T + 2k)/4$ if $k$ is even. For the $4mm_d$ 1R kolams, the minimum number of tiles needed is $(T + 4k + 3)/8$ if $k$ is odd, and $(T + 2k)/8$ if $k$ is even.

For the $2m_d m_d$ and $4mm_d$ 2R kolams, an analysis similar to that for the 1R kolams can be made, first for the $k^2$ 1R kolam and then for the $(k − 1)^2$ 1R kolam. Summing them up gives the total minimum number of tiles required to complete the kolam. The minimum number of tiles to be specified for $2mm_d$ 2R kolams is $(T+4k-1)/4$ for both $k$-odd and $k$-even. The minimum number of tiles to be specified for $4m_d m_d$ 2R kolams is $(T+6k+1)/8$ if $k$ is odd, and $(T+6k-1)/8$ if $k$ is even. 2R kolams will always have tiles lying on the mirror lines, and these tiles must have mirror symmetry about the diagonal of the individual square tiles. This restricts the tiles lying on the mirror $m_d$ in 2R kolams to either Circle, Door, or Diamond tiles.

•**1-fold rotation:** All kolams have 1-fold rotation symmetry. If this is the only symmetry present in a kolam, there are no further constraints beyond the three rules for composing the kolam.

**8 Discussion and conclusions**

We have classified all square-tile kolams into 8 symmetry groups, and presented expressions for enumerating them for two specific templates, namely, single rectangle (1R) and two rectangle (2R). We conclude by noting some future extensions of this study.

We have presented only 6 basic square tiles in Figure1; there can be many more in principle by adding flourishes to the lines that encircle the dots and by considering connectivity between dots that are not nearest neighbors. There are also square tiles with different orientations of the lines within a tile that would be needed to cover a larger variety of kolams. Assembling the square tiles in Figure 1 can only generate kolams with a square or rectangular grid of dots. In general, any arbitrary arrangement of dots should be possible. Hexagonal tiles would be a straightforward extension to this work that will be explored in a separate study. An interesting transformation between the 2R and hexagonal grid templates is that a square 2R template when stretched in one direction by √3 turns it into a regular hexagonal grid. The discussion here with the six basic square tiles is not comprehensive in terms of being able to make all types of kolams with polygonal-shaped tiles. Nonetheless, as Table 2 shows, the number of kolams that are possible with just these six tiles balloons exponentially with the number of tiles, illustrating the tremendous richness of this art form.

## 9 Acknowledgements

I would like to thank Dr. Amy Alznauer from Northwestern University who helped co-design these square tile kolams and Christopher Danielson for manufacturing these tiles from maple wood.

## 11 Appendix

| Point Group (PG) symmetry of the kolam | | | | | | | | |
|---|---|---|---|---|---|---|---|---|
| 1 | $m(\perp \mathbf{k})$ | $m(\perp \mathbf{l})$ | $m_d$ | 2 | $2mm$ | $2m_d m_d$ | 4 | $4mm_d$ |
| $E_s^{PG}$ for single rectangle (1R): $(kl)$ tiles | | | | | | | | |
| $2kl - k - l$ | $k$-even: $kl - \frac{k}{2}$ | $l$-even: $kl - \frac{l}{2}$ | $k=l=n$: $n(n-1)$ | $(k+l)$ even: $kl - \frac{k+l}{2}$ | $(kl)$ even: $\frac{kl}{2}$ | $k=l=n$: $\frac{n(n-1)}{2}$ | $k=l=n$: $\frac{n(n-1)}{2}$ | $k=l=n$ even: $\frac{n^2}{4}$<br>$k=l=n$ odd: $\frac{n^2-1}{4}$ |
| | $k$-odd: $kl - \frac{(k+1)}{2}$ | $l$-odd: $kl - \frac{(l+1)}{2}$ | $k \neq l$: no kolam | $(k+l)$ odd: $kl - \frac{(k+l-1)}{2}$ | $(kl)$ odd: $\frac{(kl-1)}{2}$ | $k \neq l$: no kolam | $k \neq l$: no kolam | $k \neq l$: no kolam |
| $E_s^{PG}$ for two rectangles (2R): $(kl + (k-1)(l-1))$ tiles | | | | | | | | |
| $4(kl - k - l + 1)$ | $2(k-1) \times (l-1)$ | $2(k-1) \times (l-1)$ | $k=l=n$: $2(n-1)^2 + (n-1)$ | $2(k-1) \times (l-1)$ | $(k-1) \times (l-1)$ | $k=l=n$: $(n-1)^2 + (n-1)$ | $k=l=n$: $(n-1)^2$ | $k=l=n$: $\frac{(n-1)^2}{2} + \frac{n-1}{2}$ |
| | | | $k \neq l$: no kolam | | | $k \neq l$: no kolam | $k \neq l$: no kolam | $k \neq l$: no kolam |

Table 1. The number of minimum shared edges, $E_s^{PG}$, for a kolam of a given template in Figure 3, classified by its symmetry group labeled by the point groups ($PG$), illustrated in Figure 2. The number of kolams is given by Equation (1). The bold labels $\mathbf{k}$ and $\mathbf{l}$ denote the respective edge directions in which $k$ tiles and $l$ tiles are counted in the rectangular templates in Figure 3.

| $k \times l$ | 1 | $m(\perp \mathbf{k})$ | $m(\perp \mathbf{l})$ | $m_d$ | 2 | $2mm$ | $2m_d m_d$ | 4 | $4mm_d$ |
|---|---|---|---|---|---|---|---|---|---|
| **Number, $N(PG)$, of single rectangle (1R) kolams for each Point Group (PG)** | | | | | | | | | |
| $1 \times 1$ | $2^0$ | $2^0$ | $2^0$ | $2^0$ | $2^0$ | $2^0$ | $2^0$ | $2^0$ | $2^0$ |
| $1 \times 2$ | $2^1$ | $2^1$ | $2^1$ | - | $2^1$ | $2^1$ | - | - | - |
| $2 \times 2$ | $2^4$ | $2^3$ | $2^3$ | $2^2$ | $2^2$ | $2^2$ | $2^1$ | $2^1$ | $2^1$ |
| $2 \times 3$ | $2^7$ | $2^5$ | $2^4$ | - | $2^4$ | $2^3$ | - | - | - |
| $3 \times 3$ | $2^{12}$ | $2^7$ | $2^7$ | $2^6$ | $2^6$ | $2^4$ | $2^3$ | $2^3$ | $2^2$ |
| $3 \times 4$ | $2^{17}$ | $2^{10}$ | $2^{10}$ | - | $2^8$ | $2^6$ | - | - | - |
| $4 \times 4$ | $2^{24}$ | $2^{14}$ | $2^{14}$ | $2^{12}$ | $2^{12}$ | $2^8$ | $2^6$ | $2^6$ | $2^4$ |
| $4 \times 5$ | $2^{31}$ | $2^{18}$ | $2^{17}$ | - | $2^{16}$ | $2^{10}$ | - | - | - |
| $5 \times 5$ | $2^{40}$ | $2^{22}$ | $2^{22}$ | $2^{20}$ | $2^{20}$ | $2^{12}$ | $2^{10}$ | $2^{10}$ | $2^6$ |
| **Number, $N(PG)$, of two rectangles (2R) kolams for each Point Group (PG)** | | | | | | | | | |
| $1 \times 1$ | $2^0$ | $2^0$ | $2^0$ | $2^0$ | $2^0$ | $2^0$ | $2^0$ | $2^0$ | $2^0$ |
| $1 \times 2$ | $2^0$ | $2^0$ | $2^0$ | - | $2^0$ | $2^0$ | - | - | - |
| $2 \times 2$ | $2^4$ | $2^2$ | $2^2$ | $2^3$ | $2^2$ | $2^1$ | $2^2$ | $2^1$ | $2^1$ |
| $2 \times 3$ | $2^8$ | $2^4$ | $2^4$ | - | $2^4$ | $2^2$ | - | - | - |
| $3 \times 3$ | $2^{12}$ | $2^8$ | $2^8$ | $2^{10}$ | $2^8$ | $2^4$ | $2^6$ | $2^4$ | $2^3$ |
| $3 \times 4$ | $2^{24}$ | $2^{12}$ | $2^{12}$ | - | $2^{12}$ | $2^6$ | - | - | - |
| $4 \times 4$ | $2^{28}$ | $2^{18}$ | $2^{18}$ | $2^{21}$ | $2^{18}$ | $2^9$ | $2^{12}$ | $2^9$ | $2^6$ |
| $4 \times 5$ | $2^{40}$ | $2^{24}$ | $2^{24}$ | - | $2^{24}$ | $2^{12}$ | - | - | - |
| $5 \times 5$ | $2^{56}$ | $2^{32}$ | $2^{32}$ | $2^{36}$ | $2^{32}$ | $2^{16}$ | $2^{20}$ | $2^{16}$ | $2^{10}$ |

Table 2. The number of unique kolams, $N(PG)$, from Equation (1) for the 1R and 2R templates depicted in Figure 3. The expressions for the minimum shared edges in Table 1 are used. The bold labels $\mathbf{k}$ and $\mathbf{l}$ denote the respective edge directions in which $k$ tiles and $l$ tiles are counted in the rectangular templates in Figure 3.